\documentclass[11pt]{article}
\usepackage{amsfonts}
\usepackage{graphicx}

\usepackage{amsmath}
\usepackage{amssymb}
\usepackage{latexsym}
\usepackage{amsmath, amsfonts,amssymb, amsthm, enumerate, euscript,makeidx,color,mathrsfs}

\def\sqr#1#2{{\vcenter{\vbox{\hrule height.#2pt
              \hbox{\vrule width.#2pt height#1pt \kern#1pt \vrule width.#2pt}
              \hrule height.#2pt}}}}

\usepackage{fancyhdr}
\usepackage{amssymb}
\usepackage{latexsym}
\usepackage{mathrsfs}
\usepackage{cancel}
\usepackage{dsfont}
\usepackage{varioref}
\usepackage{amsmath}
\usepackage{color}
\usepackage{multirow}
\usepackage{eso-pic}
\usepackage{enumerate}
\usepackage{cleveref}

\newtheorem{Definition}{Definition}[section]
\newtheorem{Theorem}[Definition]{Theorem}
\newtheorem{Lemma}[Definition]{Lemma}

\newtheorem{Proposition}[Definition]{Proposition}
\newtheorem{Remark}{Remark}[section]
\newtheorem{Example}{Example}[section]

\def\nnb{\nonumber}
\def\ds{\displaystyle}

\def\all{  \, \forall \, }

\newcommand{\refeq}[1]{$(\ref{#1})$}
\newcommand{\thb}[1]{~{\rm (#1)}}
\def\Proof{\noindent{\bf Proof. }}
\def\endpf{\hfill$\Box$\vspace{0.4cm}}
\def\eqif{\, {\rm if}\, \,}
\def\eqon{ \, {\rm on } \, \,}
\def\eqin{ \, {\rm in } \, \,}
\def\eqae{ \, {\rm a.e. } \,\, }
\newcommand{\si}[1]{\mbox{ strongly in}\, #1}
\newcommand{\wi}[1]{\mbox{ weakly  in}\, #1}

\newcommand{\ui}[1]{\mbox{uniformly in}\, #1}
\def\ol{\overline}
\def\Ga{\alpha}

\def\Gg{\gamma}
\def\Gd{\delta}

\def\Gve{\varepsilon}

\def\Gs{\sigma}

\def\Gvp{\varphi}
\def\GO{\Omega} 
\def\bY{\ol{Y}}

\def\bu{{\bar u}}
\def\bv{{\bar v}}

\def\by{{\bar y}}

\def\cA{{\cal A}}

\def\cU{{\cal U}}
\def\cV{{\cal V}}

\def\cX{{\cal X}}

\def\bpsi{{\bar \psi}}

\def\bGO{{\ol{\Omega}}}
\def\tiC{\widetilde C}

\def\tiu{\tilde u}
\def\tiv{\tilde v}

\def\tiy{\tilde y}

\def\qq{\qquad}
\def\q{\quad}
\newcommand{\pri}{{\prime}}

\def\pa{\partial}
\def\na{\nabla}
\def\Div{ \, {\rm div } \,\, } 

\def\Limsup{\varlimsup}
\def\Liminf{\varliminf}

\if{

}\fi

\def\IR{\mathbb{R}} 
\newcommand{\set}[1]{\left\{#1\right\}}
\newcommand{\ip}[1]{\left\langle #1\right\rangle}

\begin{document}
\title{Optimality conditions for optimal control of multisolution p-Laplacian elliptic equations\thanks{This work was
 supported by the National Natural Science Foundation of
China (11726619, 11726620, 11601213 and 11771097), the Natural Science Foundation of Guangdong Province
(2018A0303070012), and the Key Subject Program of Lingnan Normal University (1171518004).}
}
\author{Hongwei Lou\footnote{School
of Mathematical Sciences, and LMNS, Fudan University, Shanghai
200433, China (Email: \texttt{hwlou@fudan.edu.cn}).}
~~and~~Shu Luan\footnote{School of Mathematics and Statistics, Lingnan Normal University,
Zhanjiang, Guangdong 524048, China (Email: \texttt{luanshu@yeah.net}).}
}
\date{}

\maketitle

\begin{quote}
\footnotesize {\bf Abstract.} In this paper, an optimal control problem governed by a class of p-Laplacian elliptic equations is studied.
In particular, as no monotonicity assumption is assumed on the nonlinear term, the state equation may admit several
solutions for one control. To obtain optimality conditions for an optimal pair, the multiplicity
and singularity/degeneracy of the state equation need to be handled respectively.
For this reason, penalization problems and approximation problems  are introduced.
Finally the main result is proved by a series of process of taking to the limits.


\end{quote}

\section{Introduction}
Due to some practical interests, many authors studied optimal control for elliptic differential equations. Most of these works deal with well-defined state equations. We refer the readers to the books by Li and Yong\cite{14}, Barbu\cite{2}, Berkovitz\cite{3},
Clarke\cite{8} and the papers \cite{6, 7, 9, 11, 12, 20, 22, 24} for further details.

In the present paper, we study an optimal control problem governed by a class of non-well-defined $p$-Laplacian elliptic equations.
A review on important applications of optimal control theory to problems in engineering and medical science shows
that in most of the cases the underlying PDEs are quasilinear. The state equation (p-Laplacian elliptic equation) considered
in our paper is a typical quasilinear equation, which arises from the studies of nonlinear phenomena in non-Newtonian fluids, reaction-diffusion problems,
non-linear elasticity, torsional creep problem, glacelogy, radiation of heat, etc.(see \cite{1}). In particular, the case
where $1<p<2$ is of the most interest for elastic-plastic models (see \cite{17}). Moreover, it is pointed out that
no monotonicity assumption is posed on the nonlinear term, then the state equation may admit
 several solutions for one control and hence is non-well-posed. Such non-well-posed equations are mainly
found in bifurcation theory. Some models describing enzymatic reactions, phenomena in plasma physics and chemistry
have also this property (see Crandall and Rabinowitz \cite{10} and Lions \cite{15} for more discussions).

As we know, in the case that a state equation admits more than one solution, the state variable does not depend continuously on
the control. Therefore, we cannot obtain the variations of the state with respect to the control similarly as
in [14]. Generally, a penalization approach is considered to deal with such non-well-posed cases. In [16],
Lions first studied optimal control of non-monotone elliptic systems without state constraints, while Bonnans and Casas
[4] considered the case in which the state constraints were involved. Their methods are to penalize the problem by
removing the nonlinear term from the state equation and regarding it a part of the state constraints. In recent years,
some authors discussed  more general state equations (see \cite{5, 13, 19, 21, 23}).

The main difference between this paper and the existing literatures lies in that
the state equation we considered has both singularity/degeneracy and multiplicity, and we have to deal with
the two difficulties, respectively. For this reason we first introduce penalization problems and approximation problems.
A natural question is that if the penalization problem and approximation problem can be discussed together.
Unfortunately, it doesn't seem to work. (see Remark \ref{R301}).

The rest of this paper is organized as follows. In Section 2, we will give the formulation of the control problem and the main
result. Section 3 is devoted to constructing penalization problems and approximation problems.
In Section 4, we give optimality conditions for penalization problems and approximation problems.
Our main result will be proved in Section 5.
Finally, we give an example in Section 6 to show an application of the main result.

\section{Formulation of the control problem and the main result}
\def\theequation{2.\arabic{equation}}
\setcounter{equation}{0}
Let $1<p<2$ and $\GO$ be a bounded domain of $\IR^n (n\ge 1)$ with
$C^{1,1}$ boundary $\pa \GO$.
Denote $p^*=\frac{np}{n-p}$ if $n>p$ and $p*=+\infty$  if $n\leq p$. Or equivalently,
 $p^*=\frac{np}{n-p}$ if $n\geq 2$ and $p*=+\infty$  if $n=1$.

Consider the following $p$-Laplacian elliptic
equation
\begin{equation}\label{E201}
\left\{
\begin{array}{ll}
-{\rm div}(|\nabla y|^{p-2}\nabla y)=f(y)+u
&\eqin \GO,\\[1mm]
y=0&\eqon
\pa \GO
\end{array}
\right.
\end{equation}
and the cost functional
\begin{equation}\label{E202}
J(y,u)=\int_{\GO}\Big(f^0(x,y(x))\,dx+g(x,u(x))\Big)\, dx.
\end{equation}

We  set the following assumptions.

\vskip2mm (S1) Assume $U=[a,b]\subset\IR$. Denote by  $\cU\equiv \{u:\GO\to [a,b]\big|u \, \mbox{is measurable}\}$ the control set.

\vskip1mm (S2) The function $f\in C^1(\IR)$  satisfies the growth condition
$$
|f(y)|\le C(1+|y|^{r-1}),\qq y\in \IR,
$$
where $C>0$ and $r\in [1,p^*)$ are constants.

\vskip1mm (S3) The function $f^0: \GO\times\IR \to \IR$
satisfies the following properties: $f^0(\cdot,y)$ is measurable in $\GO$,
$f^0(x,\cdot)$ and $f^0_y(x,\cdot)$
are continuous in $\IR$, and for any $M>0$, there exists a constant $C_M>0$ such that
$$
|f^0(x,y)|+|f^0_y(x,y)|\le C_M, \qq\all x\in \GO,\, |y|\le M.
$$

\vskip1mm (S4) The function $g: \GO\times U \to \IR$
is measurable in $x\in\GO$, continuous and convex in $u\in U$. Moreover, $g$ is bounded on $\GO\times U$.

Denote
$$
\cA= \set{(y,u)\in W^{1,p}_0 (\GO)\times \cU\,
\big|(y,u)\, \mbox{satisfies \refeq{E201}}}
$$
the set of admissible pairs.

The optimal control problem is stated as follows.

\textbf{Problem (P)}. Find a pair $(\by,\bu)\in\cA$
such that
$$
J(\by,\bu)=\inf_{(y,u)\in\cA}J(y,u).
$$

A solution $(\by,\bu)\in\cA$ of Problem (P) is said to be an optimal pair,
$\bu$ is called an optimal control, and $\by$ is called an
optimal state.

The purpose of this paper is to give an optimality condition for an optimal pair
$(\by,\bu)$. \vskip2mm

\begin{Remark}\label{R201} Since no monotonicity assumption such as $f'(y)\le 0$ $(\all y\in\IR)$ is assumed,
the state
equation \refeq{E201} may admit more than one solution for some
$u\in\cU$. Hence, \refeq{E201} is non-well-posed.
\end{Remark}

By  (S2) and standard De Giorgi estimate, we can get the following proposition.
\begin{Proposition}\label{P201}  Assume that \thb{S2} holds. Then there exists a constant $C>0$, independent of $u\in \cU$,
such that $\|y\|_{L^\infty(\GO)}\le C$
for any solution $y$ of \thb{\ref{E201}}.
\end{Proposition}
For $y\in H^1_0(\GO)$, denote
\begin{equation}\label{E202B}
\cX(y):=\set{\Gvp\in H^1_0(\GO)\big| \int_{\{\na y\neq 0\}} |\na y|^{p-2}\,|\na \Gvp|^2\, dx<+\infty}.
\end{equation}
The main result of this paper is as follows.
\begin{Theorem}\label{T201} Assume that $1<p<2$ and  {\em(S1)--(S4)} hold.
 Let $(\by,\bu)\in \cA$ be an optimal
pair of Problem {\em(P)}. Then there exist a real number $\mu\ge0$ and a function $\bar{\psi}\in \cX(\by)$
such that
\begin{eqnarray}\label{E203}
\mu+\|\bar{\psi}\|_{H^1_0({\rm\Omega})}>0,
\end{eqnarray}
\begin{eqnarray}\label{E204}
& & -\Div \Big[|\nabla\by|^{p-2}\Big(I+(p-2)
\frac{\nabla\by(\nabla\by)^T}{|\nabla\by|^2}\Big)\na \bpsi\Big]= f^\pri(\by)\bpsi-\mu f^0_y(x,\by), \q\eqin \{\na \by\neq 0\}.
\end{eqnarray}
\begin{equation}\label{E205}
\na \bpsi=0,\qq\eqae \{\na \by=0\}
\end{equation}
and
\begin{equation}\label{E206}
\bpsi(x) \bu(x)-\mu g(x,\bu(x))=\max_{a\leq u\leq b}\Big(\bpsi(x) u-\mu g(x,u)\Big), \qq\eqae x\in \GO.
\end{equation}
\end{Theorem}

\begin{Remark}
Assumption (S2) is used mainly to guarantee the boundeness of $\by$ in $L^\infty(\GO)$. If we assume that $\by\in L^\infty(\GO)$, then instead of (S2), we need only to suppose that $f\in C^1(\IR)$.
\end{Remark}
\begin{Remark}
We failed to get necessary conditions for  the case of $p>2$. In {\rm\cite{Lou}}, necessary conditions for  the case of $p>2$ were only established when $f^\pri(y)\leq -\Gg$ for some constant $\Gg>0$. Yet, $f^\pri(y)\leq -\Gg$ implies that \refeq{E201} is well-posed.
\end{Remark}

\begin{Remark}
Necessary condition like Theorem \ref{T201} looks quite inadequate. Yet it still contains crucial information of the optimal pair. For example, in {\rm\cite{Lou6}}, similar result was used to analyze the regularity and existence of optimal control. While in Section 6, we give an example to show a usage of  Theorem \ref{T201}.
\end{Remark}

\section{Penalization problems and approximation problems}
\def\theequation{3.\arabic{equation}}
\setcounter{equation}{0}

To treat Problem (P), we meet
two
main difficulties. One
 is that the state equation is not well-defined. Thus, we need to
construct penalization problems
corresponding to  Problem (P) first.

Let $(\by,\bu)$ be an optimal pair of Problem (P).
Consider the following system
\begin{equation}\label{E301}
\left\{
\begin{array}{ll}
-{\rm div}(|\nabla y|^{p-2}\nabla y)=v+u
&\eqin \GO,\\[1mm]
y=0&\eqon
\pa \GO,
\end{array}
\right.
\end{equation}
where the control $(v,u)\in\cV\times\cU$ with
\begin{equation}\label{E302}
\cV\equiv\{v\in L^{\infty}(\GO):\|v-f(\by)\|_{L^{\infty}(\GO)}\le 1\}.
\end{equation}
We denoted by $y^{(v,u)}$ the solution of \refeq{E301} corresponding to $(v,u)$.

For $\tau\in(0,1)$ and $m\geq 1$, consider  the following cost functional
\begin{eqnarray}\label{E303}
\nnb J_{\tau,m}(v,u)&=& J(y^{(v,u)},u)+\int_\GO\Big(m\big|v-f(y^{(v,u)})\big|^2\\
&& +\tau|u-\bu|^2+\tau\big|v-f(\by)\big|^2\Big)\, dx.
\end{eqnarray}
We set

\bigskip
\if{
$J_{\tau,m}$ µÄ×îÓŽâ $\bv
_m$ µÄ¼«ÏÞ $\tiv$ Âú×ã $\tiv=f(\tiy)$, ¿´ÉÏÈ¥Ëƺõ·½³Ì (2.1) ×Ü»áÓнâ. ÊÂʵÉÏ, ÕâÀïºÜÖØÒªµÄÊÇ (2.1) Óнâ, ²Åµ¼ÖÂÁË $\tiv=f(\tiy)$ ³ÉÁ¢. ÕâÊÇÒòΪÓÉÓÚ (2.1) Óнâ, ʹµÃ ÏàÓ¦µÄÖ¸±ê $J_{\tau,m}$ µ± $m$ Ç÷ÓÚÎÞÇîʱÊÇÓÐÏÞµÄ.

ÕâÀïÓÐȤµÄÊÇ, (2.1) ½âµÄ´æÔÚÐÔÓë $J_{\tau,m}$ Öµº¯Êý(¹ØÓÚ $m$ µÄ)Ò»ÖÂÓнçÐԵȼÛ.
}\fi

\textbf{Problem} $(P_{\tau,m})$. Find
$(\bv_{\tau,m},\bu_{\tau,m})\in\cV\times\cU$ such that
$$
J_{\tau,m}(\bv_{\tau,m},\bu_{\tau,m})=\inf_{(v,u)\in\cV\times\cU}
J_{\tau,m}(v,u).
$$

The another main difficulty is that the state equation is singular/degenerate. Therefore,  we need to
introduce approximation problems.
For  $\Gve\in [0,1]$, consider
\begin{equation}\label{E304}
\left\{
\begin{array}{ll}
-\Div\Big(\big(\Gve^2+|\nabla y|^2\big)^{p-2\over 2}\nabla y\Big)=v +u
\qq\eqin \GO,\\[1mm]
y=0\qq\eqon \pa \GO
\end{array}
\right.
\end{equation}
and
denote by $y_\Gve^{(v,u)}$ the solution of \refeq{E304} corresponding to $(v,u)$.

Further, for an optimal control $(\bv_{\tau,m},\bu_{\tau,m})$ of Problem ($P_{\tau,m}$)  and $\Gs>0$, consider
\begin{eqnarray}\label{E305}
\nnb &&J^{\Gs,\Gve}_{\tau,m}(v,u)= J(y_\Gve^{(v,u)},u)+\int_\GO\Big(m\big|v-f(y_\Gve^{(v,u)})\big|^2+\tau\big|v-f(\by)\big|^2\\
&&\qq+\tau|u-\bu|^2+\Gs |v-\bv_{\tau,m}|^2+\Gs |u-\bu_{\tau,m}|^2\Big)\, dx
\end{eqnarray}
and

 \bigskip

\textbf{Problem} $(P^{\Gs,\Gve}_{\tau,m})$. Find
$(\bv^{\Gs,\Gve}_{\tau,m},\bu^{\Gs,\Gve}_{\tau,m})\in\cV\times\cU$ such that
$$
J^{\Gs,\Gve}_{\tau,m}(\bv^{\Gs,\Gve}_{\tau,m},\bu^{\Gs,\Gve}_{\tau,m})=\inf_{(v,u)\in\cV\times\cU}
J^{\Gs,\Gve}_{\tau,m}(v,u).
$$

\begin{Remark}\label{R301}
It is natural to ask if we can treat the two difficulties simultaneously. For example we consider simply Problem $(P^{\tau,{1\over m}}_{\tau,m})$ directly. The pity is that it does not seem to work. The reason is mainly that we do not know if  optimal controls of  Problem $(P^{\tau,{1\over m}}_{\tau,m})$ converge to $(f(\by),\bu)$.
\end{Remark}

By Proposition \ref{P201}, we have $f(\by)\in L^\infty(\GO)$. Thus, the following lemma becomes a special case of Theorem 1 in \cite{Lieberman}, which shows the
existence, uniqueness and regularity of the solution for \refeq{E304} (especially, \refeq{E301} when $\Gve=0$).
\begin{Lemma}\label{L301}
Assume that {\em (S1)} and {\em (S2)} hold. Then for any $\Gve \in [0,1],
~(v,u)\in\cV\times\cU$, \refeq{E304} admits a unique solution
$y_\Gve^{(v,u)}\in W^{1,p}_0(\GO)$.  Moreover, there exist constants $C>0$ and $\alpha\in (0,1)$,
independent of $\Gve\in [0,1]$ and $(v,u)\in\cV\times\cU$, such that
\begin{equation}\label{E306}
\|y_\Gve^{(v,u)}\|_{C^{1,\alpha}(\bGO)}\le C.
\end{equation}
\end{Lemma}

\bigskip

The following lemma shows the existence of an optimal control for Problem $(P_{\tau,m})$.
\begin{Lemma}\label{L302}
Assume that \thb{S1}--{\rm \,(S4)} hold. Then
Problem $(P_{\tau,m})$ admits  at least one optimal control
$(\bv_{\tau,m},\bu_{\tau,m})\in\cV\times\cU$.
\end{Lemma}

\Proof By Lemma \ref{L301}, there exist constants $C>0$ and $\alpha\in(0,1)$,
 such that for any $(v,u)\in\cV\times\cU$,
\begin{equation}\label{E307}
\|y^{(v,u)}\|_{C^{1,\alpha}(\bGO)}\le C.
\end{equation}
Thus, it follows from (S3) that
\begin{equation}\label{E308}
 \inf_{(v,u)\in\cV\times \cU}J_{\tau,m}(v,u)>-\infty.
\end{equation}
Hence, there exists a minimizing sequence
$(u_{\tau,m,k},v_{\tau,m,k})\in\cV\times\cU$ such that:
\begin{equation}\label{E309}
\lim_{k\to +\infty}J_{\tau,m}(v_{\tau,m,k},u_{\tau,m,k})
=\inf_{(v,u)\in\cV\times\cU}
J_{\tau,m}(v,u).
\end{equation}
Denote $y_{\tau,m,k}=y^{(v_{\tau,m,k},u_{\tau,m,k})}$. Then by \refeq{E307} and Arzel\'a-Ascoli's theorem,
we have that along a subsequence of  $k\to  +\infty$,
\begin{equation}\label{E310}
y_{\tau,m,k}\to \by_{\tau,m}, \qq \ui{C^1(\bGO)}.
\end{equation}
Moreover, by (S1) and the definitions of $\cV$ and $\cU$,  we have that $v_{\tau,m,k}$ and $u_{\tau,m,k}$ are bounded uniformly in $L^2(\GO)$ with respect to $\tau, m$ and $k$.
Thus, along a subsequence of $k\to +\infty$, we have
\begin{equation}\label{E311}
v_{\tau,m,k}\to  \bv_{\tau,m},\q
u_{\tau,m,k}\to \bu_{\tau,m},\q \wi{L^2(\GO)}.
\end{equation}
In addition, it is easy to see that
$\bv_{\tau,m}\in\cV$ and $\bu_{\tau,m}\in\cU$ since $\cV$ and $\cU$ are convex and closed in $L^2(\GO)$. Finally, by \refeq{E310}  and \refeq{E311}, we can deduce easily that
$\by_{\tau,m}=y^{(\bv_{\tau,m},\bu_{\tau,m})}$.

On the other hand, bu Mazur's Theorem, there exist
$N_k\geq 1$ and $\set{\Ga_{k,j}\big| 1\leq j\leq N_k}$ $(\all k\geq 1)$ such that $\Ga_{k,j}\geq 0$, $\ds \sum^{N_k}_{j=1}\Ga_{k,j}=1$ $(\all k\geq 1, \, 1\leq j\leq N_k)$ and
$$
\tiu_{\tau,m,k}\to \bu_{\tau,m},\q \si{L^2(\GO)}.
$$
Consequently,
\begin{eqnarray*}
\nnb && \int_\GO g(x,\bu_{\tau,m})\, dx=\lim_{k\to +\infty}\int_\GO g(x,\tiu_{\tau,m,k})\, dx\\
&\leq & \Limsup_{k\to +\infty}\sum^{N_k}_{j=1}\Ga_{k,j}\int_\GO g(x,u_{\tau,m,k+j})\, dx\\
&\leq & \Limsup_{k\to +\infty}\int_\GO g(x,u_{\tau,m,k})\, dx.
\end{eqnarray*}
Moreover, replacing $u_{\tau,m,k}$ by a subsequence of it, we can get
\begin{equation}\label{E312A}
\int_\GO g(x,\bu_{\tau,m})\, dx\leq  \Liminf_{k\to +\infty}\int_\GO g(x,u_{\tau,m,k})\, dx.
\end{equation}
Thus,
\begin{eqnarray*}
 J_{\tau,m}(\bv_{\tau,m},\bu_{\tau,m})&\leq & \lim_{k\to +\infty}J_{\tau,m}(v_{\tau,m,k},u_{\tau,m,k})
=\inf_{(v,u)\in\cV\times\cU}
J_{\tau,m}(v,u)
\end{eqnarray*}
since
$$
\int_{\GO}f^0(x,\by_{\tau,m})\,dx=\lim_{k\to +\infty}\int_{\GO}f^0(x,y_{\tau,m,k})\,dx
$$
and similar to \refeq{E312A},
\begin{eqnarray*}
&&\int_{\GO}\Big(g(
x,\bu_{\tau,m})+m\big|\bv_{\tau,m}-f(\by_{\tau,m})\big|^2+\tau\big|\bv_{\tau,m}-f(\by)\big|^2+\tau|\bu_{\tau,m}-\bu|^2\Big)\, dx\\
&\leq & \Liminf_{k\to +\infty}\int_{\GO}\Big(g(x,u_{\tau,m,k})+m\big|v_{\tau,m,k}-f(y_{\tau,m,k})\big|^2+\tau\big|v_{\tau,m,k}-f(\by)\big|^2  +\tau |u_{\tau,m,k}-\bu|^2\Big)\, dx.
\end{eqnarray*}
Therefore,
$(\bv_{\tau,m},\bu_{\tau,m})$ is an optimal control for Problem $(P_{\tau,m})$.
\endpf

Similarly, we have
\begin{Lemma}\label{L303}
Assume that \thb{S1}--{\rm \,(S3)} hold. Then for any
$m\geq 1$, $\Gs>0$ and $\tau,\,\Gve\in (0,1)$,
Problem $(P^{\Gs,\Gve}_{\tau,m})$ admits  at least one optimal control
$(\bv^{\Gs,\Gve}_{\tau,m},\bu^{\Gs,\Gve}_{\tau,m})\in\cV\times\cU$.
\end{Lemma}
\section {Optimality conditions for penalization problems and approximation problems}
\def\theequation{4.\arabic{equation}}
\setcounter{equation}{0}

In this section, we will give the optimality conditions for the optimal control of   Problems
$(P_{\tau,m})$ and $(P^{\Gs,\Gve}_{\tau,m})$. Some results can be looked as special cases of those in \cite{Lou}. Nevertheless, for readers' convenience,
we will give the structures of  the proofs for these results.

We first state the result for Problem $(P^{\Gs,\Gve}_{\tau,m})$.
\begin{Proposition}\label{P401} Assume that {\rm (S1)\,--\,(S3)} hold. Let
$(\bv^{\Gs,\Gve}_{\tau,m},\bu^{\Gs,\Gve}_{\tau,m})$ be an optimal control of Problem
$(P^{\Gs,\Gve}_{\tau,m})$ and $\by^{\Gs,\Gve}_{\tau,m}$ be the corresponding optimal state. Then, there exists
a function $\bar{\psi}^{\Gs,\Gve}_{\tau,m}\in H^1_0(\GO)$ satisfying
\begin{equation}\label{E401A}
\left\{
\begin{array}{ll}
-\Div \left[\Big(\Gve^2+|\nabla\by^{\Gs,\Gve}_{\tau,m}|^2\Big)^{p-2\over 2}\left(I+(p-2)
\frac{\nabla\by^{\Gs,\Gve}_{\tau,m}(\nabla\by^{\Gs,\Gve}_{\tau,m})^T}
{\Gve^2+|\nabla\by^{\Gs,\Gve}_{\tau,m}|^2}\right)\nabla\bpsi^{\Gs,\Gve}_{\tau,m}\right]\\[3mm]
=-f^0_y(x,\by^{\Gs,\Gve}_{\tau,m})
-2m \big(f(\by^{\Gs,\Gve}_{\tau,m})-\bv^{\Gs,\Gve}_{\tau,m}\big)f^\pri(\by^{\Gs,\Gve}_{\tau,m})\qq\eqin \GO,\\[3mm]
\bpsi^{\Gs,\Gve}_{\tau,m}=0\qq\eqon \pa \GO
\end{array}
\right.
\end{equation}
such that
\begin{eqnarray}\label{E402A}
\nnb  &  &  \int_\GO \Big(H^{\Gs,\Gve}_{\tau,m}\big(x,\by^{\Gs,\Gve}_{\tau,m},\bpsi^{\Gs,\Gve}_{\tau,m}, \bv^{\Gs,\Gve}_{\tau,m},\bu^{\Gs,\Gve}_{\tau,m}\big)-
H^{\Gs,\Gve}_{\tau,m}\big(x,\by^{\Gs,\Gve}_{\tau,m},\bpsi^{\Gs,\Gve}_{\tau,m}, v, u\big)\Big)\, dx\geq 0,\\
&& \qq\qq\qq\qq\qq\all (v,u)\in \cV\times \cU.
\end{eqnarray}
where
\begin{eqnarray}\label{E403A}
\nnb H^{\Gs,\Gve}_{\tau,m}(x,y,\psi, v, u)&=& \psi (v+u)-g(x,u)-m|v-f(y)|^2-\tau|u-\bu(x)|^2\\
\nnb &&  -\tau|v-f(\by(x))|^2-\Gs|u-\bu_{\tau,m}(x)|^2-\Gs|v-\bv_{\tau,m}(x)|^2, \\
&&  \q \all (x,y,\psi, v, u)\in \IR^n\times \IR\times \IR\times \IR\times \IR.
\end{eqnarray}
\end{Proposition}
\Proof
Let $(v,u)\in\cV\times\cU$. For $\delta\in (0,1)$ and $k\geq 1$, we set
$$
(v^k_\Gd(x),u^k_\Gd(x))
:=\left\{\begin{array}{ll} (v(x), u(x)), & \eqif \set{kx_1}\in [0,\Gd), \\
(\bv^{\Gs,\Gve}_{\tau,m}(x), \bu^{\Gs,\Gve}_{\tau,m}(x)), & \eqif \set{kx_1}\in [\Gd,1),
\end{array}\right. \qq\eqin \GO,
$$
where $x=(x_1,x_2,\ldots,x_n)$ and $\set{a}$ denote the decimal part of a real number $a$. Then $(v^k_\Gd,u^k_\Gd)\in \cV\times\cU$ and
it is not difficult to see that as $k\to +\infty$,
$$
y_\Gve^{(v^k_\Gd,u^k_\Gd)}\to y_\Gd\qq \ui{C^1(\bGO)}, \wi{W^{1,p}_0(\GO)}
$$
with
\begin{equation}\label{E404}
\left\{
\begin{array}{l}
-\Div\left(\Big(\Gve^2+|\nabla y_\Gd|^2\Big)^{p-2\over 2}\nabla y_\Gd\right)\\
\qq\qq\qq =(1-\Gd) \big(\bv^{\Gs,\Gve}_{\tau,m}+\bu^{\Gs,\Gve}_{\tau,m}\big)+\Gd (v+u),\q \eqin \GO,\\
y_\Gd=0 \qq \eqon \pa \GO.
\end{array}
\right.
\end{equation}
Furthermore, we have
\begin{eqnarray}
\nnb && J^{\Gs,\Gve}_{\tau,m}(\bv^{\Gs,\Gve}_{\tau,m},\bu^{\Gs,\Gve}_{\tau,m})\leq J^\Gd\equiv \lim_{k\to +\infty}J^{\Gs,\Gve}_{\tau,m}(v^k_\Gd,u^k_\Gd)\\
\nnb &=&(1-\Gd)  \int_\GO \Big(f^0(x,y_\Gd)+g(x,\bu^{\Gs,\Gve}_{\tau,m})+m\big|\bv^{\Gs,\Gve}_{\tau,m}-f(y_\Gd)\big|^2+\tau\left|\bv^{\Gs,\Gve}_{\tau,m}-f(\by)\right|^2\\
\nnb &&\qq +\tau |\bu^{\Gs,\Gve}_{\tau,m}-\bu|^2+\Gs|\bv^{\Gs,\Gve}_{\tau,m}-\bv_{\tau,m}|^2+\Gs |\bu^{\Gs,\Gve}_{\tau,m}-\bu_{\tau,m}|^2\Big)\, dx\\
\nnb && +\Gd  \int_\GO \Big(f^0(x,y_\Gd)+g(x,u)+m\big|v-f(y_\Gd)\big|^2+\tau\left|v-f(\by)\right|^2\\
\nnb &&\qq +\tau |u-\bu|^2+\Gs|v-\bv_{\tau,m}|^2+\Gs |u-\bu_{\tau,m}|^2\Big)\, dx.
\end{eqnarray}
Therefore,
\begin{eqnarray}\label{E405}
\nnb && 0\leq  \lim_{\Gd\to 0^+}{J^\Gd-J^{\Gs,\Gve}_{\tau,m}(\bv^{\Gs,\Gve}_{\tau,m},\bu^{\Gs,\Gve}_{\tau,m})\over \Gd}\\
\nnb &=&\lim_{\Gd\to 0^+}\int_\GO {f^0(x,y_\Gd)-f^0(x,\by^{\Gs,\Gve}_{\tau,m})\over \Gd}\, dx\\
\nnb && +\lim_{\Gd\to 0^+}(1-\Gd)m\int_\GO {\big|\bv^{\Gs,\Gve}_{\tau,m}-f(y_\Gd)\big|^2-\big|\bv^{\Gs,\Gve}_{\tau,m}-f(\by^{\Gs,\Gve}_{\tau,m})\big|^2\over \Gd}\, dx\\
\nnb &&+ \lim_{\Gd\to 0^+}\int_\GO\Big[\Big(g(x,u)+m\big|v-f(y_\Gd)\big|^2+\tau\big|v-f(\by)\big|^2+\tau |u-\bu|^2+\Gs|v-\bv_{\tau,m}|^2\\
\nnb &&\qq\qq\q+\Gs |u-\bu_{\tau,m}|^2\Big)-\Big(g(x,\bu^{\Gs,\Gve}_{\tau,m})+m\big|\bv^{\Gs,\Gve}_{\tau,m}-f(\by^{\Gs,\Gve}_{\tau,m})\big|^2
+\tau\big|\bv^{\Gs,\Gve}_{\tau,m}-f(\by)\big|^2\\
\nnb &&\qq\qq\q +\tau|\bu^{\Gs,\Gve}_{\tau,m}-\bu|^2+\Gs|\bv^{\Gs,\Gve}_{\tau,m}-\bv_{\tau,m}|^2+\Gs |\bu^{\Gs,\Gve}_{\tau,m}-\bu_{\tau,m}|^2\Big)\Big]\, dx
 \\
\nnb &=&\int_\GO f^0_y(x,\by^{\Gs,\Gve}_{\tau,m}) \bY^{\Gs,\Gve}_{\tau,m}\, dx+2m\int_\GO \big(f(\by^{\Gs,\Gve}_{\tau,m})-\bv^{\Gs,\Gve}_{\tau,m}\big)f^\pri(\by^{\Gs,\Gve}_{\tau,m}) \bY^{\Gs,\Gve}_{\tau,m}\, dx\\
\nnb && +\int_\GO\Big[\Big(g(x,u)+m\big|v-f(\by^{\Gs,\Gve}_{\tau,m})\big|^2+\tau\big|v-f(\by)\big|^2+\tau |u-\bu|^2+\Gs|v-\bv_{\tau,m}|^2\\
\nnb &&\qq\qq\q+\Gs |u-\bu_{\tau,m}|^2\Big)-\Big(g(x,\bu^{\Gs,\Gve}_{\tau,m})+m\big|\bv^{\Gs,\Gve}_{\tau,m}-f(\by^{\Gs,\Gve}_{\tau,m})\big|^2+\tau\big|\bv^{\Gs,\Gve}_{\tau,m}-f(\by)\big|^2\\
\nnb &&\qq\qq\q +\tau|\bu^{\Gs,\Gve}_{\tau,m}-\bu|^2+\Gs|\bv^{\Gs,\Gve}_{\tau,m}-\bv_{\tau,m}|^2+\Gs |\bu^{\Gs,\Gve}_{\tau,m}-\bu_{\tau,m}|^2\Big)\Big]\, dx\\
&&
\end{eqnarray}
with $\ds\bY^{\Gs,\Gve}_{\tau,m}$, which is the limit of $\ds {y_\Gd-\by^{\Gs,\Gve}_{\tau,m}\over \Gd}$ in $H^1_0(\GO)$, being the solution of the following equation:
\begin{equation}\label{E406}
\left\{
\begin{array}{l}
-\Div\left[\Big(\Gve^2+|\nabla\by^{\Gs,\Gve}_{\tau,m}|^2\Big)^{p-2\over 2}\Big(I+(p-2)
 {\nabla \by^{\Gs,\Gve}_{\tau,m}(\nabla\by^{\Gs,\Gve}_{\tau,m})^T\over \Gve^2+|\nabla\by^{\Gs,\Gve}_{\tau,m}|^2}\Big)\na\bY^{\Gs,\Gve}_{\tau,m}\right]\\[4mm]
\qq\qq =v+u-\bv^{\Gs,\Gve}_{\tau,m}-\bu^{\Gs,\Gve}_{\tau,m}\qq \eqin \GO,\\[4mm]
\bY^{\Gs,\Gve}_{\tau,m}=0\qq \eqon \pa \GO.
\end{array}\right.
\end{equation}

Let $\bpsi^{\Gs,\Gve}_{\tau,m}$ be the solution of \refeq{E401A},  then it follows from \refeq{E405} that
\begin{eqnarray}\label{E407}
\nnb 0&\leq & \int_\GO \bpsi^{\Gs,\Gve}_{\tau,m}\big(\bv^{\Gs,\Gve}_{\tau,m}+\bu^{\Gs,\Gve}_{\tau,m}-v-u\big)\, dx\\
\nnb && +\int_\GO\Big[\Big(g(x,u)+m\big|v-f(\by^{\Gs,\Gve}_{\tau,m})\big|^2+\tau\big|v-f(\by)\big|^2+\tau |u-\bu|^2+\Gs|v-\bv_{\tau,m}|^2\\
\nnb &&\qq\qq\q+\Gs |u-\bu_{\tau,m}|^2\Big)-\Big(g(x,\bu^{\Gs,\Gve}_{\tau,m})+m\big|\bv^{\Gs,\Gve}_{\tau,m}-f(\by^{\Gs,\Gve}_{\tau,m})\big|^2+\tau\big|\bv^{\Gs,\Gve}_{\tau,m}-f(\by)\big|^2\\
\nnb &&\qq\qq\q +\tau|\bu^{\Gs,\Gve}_{\tau,m}-\bu|^2+\Gs|\bv^{\Gs,\Gve}_{\tau,m}-\bv_{\tau,m}|^2+\Gs |\bu^{\Gs,\Gve}_{\tau,m}-\bu_{\tau,m}|^2\Big)\Big]\, dx\\
\nnb &=& \int_\GO \Big(H^{\Gs,\Gve}_{\tau,m}\big(x,\by^{\Gs,\Gve}_{\tau,m},\bpsi^{\Gs,\Gve}_{\tau,m}, \bv^{\Gs,\Gve}_{\tau,m},\bu^{\Gs,\Gve}_{\tau,m}\big)-
H^{\Gs,\Gve}_{\tau,m}\big(x,\by^{\Gs,\Gve}_{\tau,m},\bpsi^{\Gs,\Gve}_{\tau,m}, v, u\big)\Big)\, dx,\\
&& \qq\qq\qq\qq\qq\all (v,u)\in \cV\times \cU.
\end{eqnarray}
We get the proof. \endpf

The following proposition shows that the optimal control for Problem $(P^{\Gs,\Gve}_{\tau,m})$  converges to that for  Problem $(P_{\tau,m})$.

\begin{Proposition}\label{P402} Assume that {\em (S1)--(S3)} hold. Let
$(\bv^{\Gs,\Gve}_{\tau,m},\bu^{\Gs,\Gve}_{\tau,m})$ be an optimal control of Problem
$(P^{\Gs,\Gve}_{\tau,m})$ and $\by^{\Gs,\Gve}_{\tau,m}$ be the corresponding optimal state. Then, it holds that
as $\Gve\to 0^+$,
\begin{equation}\label{E408}
\left\{\begin{array}{l}
\ds \by^{\Gs,\Gve}_{\tau,m}\to \by_{\tau,m}, \qq\ui{\,C^1(\bGO)}, \\
\ds \bv^{\Gs,\Gve}_{\tau,m}\to \bv_{\tau,m}, \qq\si{\,L^2(\GO)}, \\
\ds \bu^{\Gs,\Gve}_{\tau,m}\to \bu_{\tau,m}, \qq\si{\,L^2(\GO)}.
\end{array}\right.
\end{equation}
\end{Proposition}
\Proof
It follows from the definition of $\cV$ and $\cU$ that $\bv^{\Gs,\Gve}_{\tau,m}$ and $\bu^{\Gs,\Gve}_{\tau,m}$ are uniformly bounded in $L^2(\GO)$ respect to $\Gve\in (0,1)$. On the other hand, by Lemma \ref{L301},
 $\by^{\Gs,\Gve}_{\tau,m}$ are uniformly bounded in $C^{1,\Ga}(\bGO)$. Thus, it suffices to prove that \refeq{E408} holds along a subsequence of $\Gve\to 0^+$.

Along a subsequence of $\Gve\to 0^+$, it holds that
\begin{equation}\label{E409}
\left\{\begin{array}{l}
\ds \by^{\Gs,\Gve}_{\tau,m}\to \tiy_{\tau,m}, \qq\ui{\,C^1(\bGO)}, \\
\ds \bv^{\Gs,\Gve}_{\tau,m}\to \tiv_{\tau,m}, \qq\wi{\,L^2(\bGO)}, \\
\ds \bu^{\Gs,\Gve}_{\tau,m}\to \tiu_{\tau,m}, \qq\wi{\,L^2(\bGO)}.
\end{array}\right.
\end{equation}
It is easy to see that $\tiy_{\tau,m}=y^{(\tiv_{\tau,m},\tiu_{\tau,m})}$.

By the optimality of $(\bv^{\Gs,\Gve}_{\tau,m},\bu^{\Gs,\Gve}_{\tau,m})$,
\begin{equation}\label{E410}
J^{\Gs,\Gve}_{\tau,m}(\bv^{\Gs,\Gve}_{\tau,m},\bu^{\Gs,\Gve}_{\tau,m})\leq J^{\Gs,\Gve}_{\tau,m}(\bv_{\tau,m},\bu_{\tau,m}).
\end{equation}
By \refeq{E409}--\refeq{E410},
\begin{eqnarray}\label{E411}
\nnb && J_{\tau,m}(\tiv_{\tau,m},\tiu_{\tau,m})+\Gs \Limsup_{\Gve\to 0^+}\int_\GO \Big(|\bv^{\Gs,\Gve}_{\tau,m}-\bv_{\tau,m}|^2+|\bu^{\Gs,\Gve}_{\tau,m}-\bu_{\tau,m}|^2\Big)\, dx\\
\nnb &=& \int_\GO \Big(f^0(x,\tiy_{\tau,m})+g(x,\tiu_{\tau,m})+m\left|\tiv_{\tau,m}-f(\tiy_{\tau,m})\right|^2+\tau|\tiu_{\tau,m}-\bu|^2+\tau\big|\tiv_{\tau,m}-f(\by)\big|^2\Big)\, dx\\
\nnb &&+\Gs\Limsup_{\Gve\to 0^+}\int_\GO \Big(|\bv^{\Gs,\Gve}_{\tau,m}-\bv_{\tau,m}|^2+|\bu^{\Gs,\Gve}_{\tau,m}-\bu_{\tau,m}|^2\Big)\, dx\\
\nnb &\leq & \lim_{\Gve\to 0^+}\int_\GO f^0(x,\by^{\Gs,\Gve}_{\tau,m})\, dx\\
\nnb && +\Liminf_{\Gve\to 0^+}\int_\GO \Big(g(x,\bu^{\Gs,\Gve}_{\tau,m})+m\left|\bv^{\Gs,\Gve}_{\tau,m}-f(\by^{\Gs,\Gve}_{\tau,m})\right|^2 +\tau\big|\bv^{\Gs,\Gve}_{\tau,m}-f(\by)\big|^2+\tau|\bu^{\Gs,\Gve}_{\tau,m}-\bu|^2\Big)\, dx\\
\nnb &&+\Gs\Limsup_{\Gve\to 0^+}\int_\GO \Big(|\bv^{\Gs,\Gve}_{\tau,m}-\bv_{\tau,m}|^2+|\bu^{\Gs,\Gve}_{\tau,m}-\bu_{\tau,m}|^2\Big)\, dx\\
\nnb &\leq & \Limsup_{\Gve\to 0^+}J^{\Gs,\Gve}_{\tau,m}(\bv^{\Gs,\Gve}_{\tau,m},\bu^{\Gs,\Gve}_{\tau,m})\leq \lim_{\Gve\to 0^+}J^{\Gs,\Gve}_{\tau,m}(\bv_{\tau,m},\bu_{\tau,m})\\
&=&J_{\tau,m}(\bv_{\tau,m},\bu_{\tau,m})\leq J_{\tau,m}(\tiv_{\tau,m},\tiu_{\tau,m}).
\end{eqnarray}
This implies
$$
\Limsup_{\Gve\to 0^+}\int_\GO \Big(|\bv^{\Gs,\Gve}_{\tau,m}-\bv_{\tau,m}|^2+|\bu^{\Gs,\Gve}_{\tau,m}-\bu_{\tau,m}|^2\Big)\, dx=0.
$$
That is (as $\Gve\to 0^+$),
$$
\left\{\begin{array}{l}
\ds \bv^{\Gs,\Gve}_{\tau,m}\to \bv_{\tau,m}, \qq\si{\,L^2(\bGO)}, \\
\ds \bu^{\Gs,\Gve}_{\tau,m}\to \bu_{\tau,m}, \qq\si{\,L^2(\bGO)}.
\end{array}\right.
$$
Consequently, $\tiy_{\tau,m}=\by_{\tau,m}$ and \refeq{E408} holds.
\endpf

Now, we state the necessary  conditions  for optimal control of Problem $(P_{\tau,m})$. We have
\begin{Proposition}\label{P403} Assume that {\em (S1)--(S3)} hold and $1<p<2$. Let
$(\bv_{\tau,m},\bu_{\tau,m})$ be an optimal control of Problem
$(P_{\tau,m})$ and $\by_{\tau,m}$ be the corresponding optimal state. Then, there exists
a function $\bpsi_{\tau,m}\in H^1_0(\GO)$ satisfying
\begin{eqnarray}\label{E412C}
\nnb &&\|\bpsi_{\tau,m}\|^2_{H^1_0(\GO)}+\int_{\set{\nabla\by_{\tau,m}\ne 0}}|\nabla\by_{\tau,m}|^{p-2}\, |\nabla\bpsi_{\tau,m}|^2\, dx\\
&\leq & C\Big(1+\Big|\int_\GO 2m \big(f(\by_{\tau,m})-\bv_{\tau,m}\big)f^\pri(\by_{\tau,m})\bpsi_{\tau,m}\, dx\Big|\Big),
\end{eqnarray}
\begin{eqnarray}\label{E413A}
\nnb && \ds
-\Div\Big[|\nabla\by_{\tau,m}|^{p-2}\Big(I+(p-2)
\frac{\nabla\by_{\tau,m}(\nabla\by_{\tau,m})^T}
{|\nabla\by_{\tau,m}|^2}\Big)\nabla\bpsi_{\tau,m}\Big]\\
\nnb &=& -f^0_y(x,\by_{\tau,m})
-2m \big(f(\by_{\tau,m})-\bv_{\tau,m}\big)f^\pri(\by_{\tau,m})\\
&&\qq\qq\qq\qq\qq\qq\qq \eqin \set{\na\by_{\tau,m}\ne 0},\\
\nnb &&\\
\label{E414A} && \na \bpsi_{\tau,m}=0\qq\qq\qq\qq\qq\eqae \set{\na\by_{\tau,m}=0}
\end{eqnarray}
and
\begin{eqnarray}\label{E415A}
\nnb  &  &  \int_\GO \Big(H_{\tau,m}\big(x,\by_{\tau,m},\bpsi_{\tau,m}, \bv_{\tau,m}, \bu_{\tau,m}\big)-
H_{\tau,m}\big(x,\by_{\tau,m},\bpsi_{\tau,m}, v, u\big)\Big)\, dx\geq 0\\
&& \qq\qq\qq\qq\qq\all (v,u)\in \cV\times \cU.
\end{eqnarray}
where $C>0$ is a constant independent of $\tau>0$ and $m\geq 1$,
\begin{eqnarray}\label{E416A}
\nnb &&  H_{\tau,m}(x,y,\psi, v, u)\\
\nnb &=& \psi (v+u)-g(x,u)-m|v-f(y)|^2-\tau |v-f(\by(x))|^2-\tau |u-\bu(x)|^2, \\
&&  \qq\qq\qq \all (x,y,\psi, v, u)\in \IR^n\times \IR\times \IR\times \IR\times \IR.
\end{eqnarray}
\end{Proposition}

\Proof By \refeq{E401A}, since $1<p<2$ and $\by^{\Gs,\Gve}_{\tau,m}$ is bounded uniformly in $C^{1,\Ga}(\bGO)$, there exists  constants $\Gg>0$ and $C_1>0$, independent of $\Gs>0, m\geq 1$ and $\tau,\,\Gve\in (0,1)$,  such that
\begin{eqnarray}\label{E417B}
\nnb && \Gg \|\bpsi^{\Gs,\Gve}_{\tau,m}\|^2_{H^1_0(\GO)}\leq  (p-1)\int_\GO  \Big(\Gve^2+|\nabla\by^{\Gs,\Gve}_{\tau,m}|^2\Big)^{p-2\over 2}\big|\nabla\bpsi^{\Gs,\Gve}_{\tau,m}\big|^2\, dx\\
\nnb &\leq & \int_\GO \ip{\Big(\Gve^2+|\nabla\by^{\Gs,\Gve}_{\tau,m}|^2\Big)^{p-2\over 2}\left(I+(p-2)
\frac{\nabla\by^{\Gs,\Gve}_{\tau,m}(\nabla\by^{\Gs,\Gve}_{\tau,m})^T}
{\Gve^2+|\nabla\by^{\Gs,\Gve}_{\tau,m}|^2}\right)\nabla\bpsi^{\Gs,\Gve}_{\tau,m},\na \bpsi^{\Gs,\Gve}_{\tau,m}}\, dx\\
\nnb &=& -\int_\GO \Big(f^0_y(x,\by^{\Gs,\Gve}_{\tau,m})
+2m \big(f(\by^{\Gs,\Gve}_{\tau,m})-\bv^{\Gs,\Gve}_{\tau,m}\big)f^\pri(\by^{\Gs,\Gve}_{\tau,m})\Big)\bpsi^{\Gs,\Gve}_{\tau,m}\, dx\\
&\leq &  C_1\Big(\|\bpsi^{\Gs,\Gve}_{\tau,m}\|_{L^1(\GO)}+\Big|\int_\GO 2m \big(f(\by^{\Gs,\Gve}_{\tau,m})-\bv^{\Gs,\Gve}_{\tau,m}\big)f^\pri(\by^{\Gs,\Gve}_{\tau,m})\bpsi^{\Gs,\Gve}_{\tau,m}\, dx\Big|\Big).
\end{eqnarray}
Then, combining the above with Poincar\'e's inequality, we have
\begin{equation}\label{E418}
\|\bpsi^{\Gs,\Gve}_{\tau,m}\|^2_{H^1_0(\GO)}\leq C\Big(1+\Big|\int_\GO  2m \big(f(\by^{\Gs,\Gve}_{\tau,m})-\bv^{\Gs,\Gve}_{\tau,m}\big)f^\pri(\by^{\Gs,\Gve}_{\tau,m})\bpsi^{\Gs,\Gve}_{\tau,m}\, dx\Big|\Big)
\end{equation}
 for  some constant  $C>0$  independent of  $\Gs>0, m\geq 1$ and $\tau,\,\Gve\in (0,1)$.
And consequently, $\bpsi^{\Gs,\Gve}_{\tau,m}$ is bounded uniformly in $H^1_0(\GO)$ respect to $\Gs>0$ and $\Gve\in (0,1)$. Thus, as $\Gve\to 0^+$, we can suppose that
$\bpsi^{\Gs,\Gve}_{\tau,m}$ converges to some $\bpsi^\Gs_{\tau,m}$ weakly in $H^1_0(\GO)$ and strongly in $L^2(\GO)$.  Then, it follows from \refeq{E408} and \refeq{E418}
that
\begin{eqnarray}\label{E419}
\nnb &&\|\bpsi^\Gs_{\tau,m}\|^2_{H^1_0(\GO)}+\int_{\set{\nabla\by_{\tau,m}\ne 0}}|\nabla\by_{\tau,m}|^{p-2}\, |\nabla\bpsi^\Gs_{\tau,m}|^2\, dx\\
&\leq & C\Big(1+\Big|\int_\GO 2m \big(f(\by_{\tau,m})-\bv_{\tau,m}\big)f^\pri(\by_{\tau,m})\bpsi^\Gs_{\tau,m}\, dx\Big|\Big).
\end{eqnarray}
 Moreover, by \refeq{E408}  and the fact of that
$\set{\na \by_{\tau,m}\ne 0}:=\set{x\in\GO\big| \na\by_{\tau,m}(x)\ne 0}$ is an open subset of $\GO$, it is not difficult to get
\begin{eqnarray}\label{E420}
\nnb && \ds
-\Div\Big[|\nabla\by_{\tau,m}|^{p-2}\Big(I+(p-2)
\frac{\nabla\by_{\tau,m}(\nabla\by_{\tau,m})^T}
{|\nabla\by_{\tau,m}|^2}\Big)\nabla\bpsi^\Gs_{\tau,m}\Big]\\
 &=& -f^0_y(x,\by_{\tau,m})
-2m \big(f(\by_{\tau,m})-\bv_{\tau,m}\big)f^\pri(\by_{\tau,m})\qq \eqin \set{\na\by_{\tau,m}\ne 0}
\end{eqnarray}
from \refeq{E401A}.

By \refeq{E402A} (see also \refeq{E407}), we have
\begin{eqnarray}\label{E421}
\nnb 0&\leq & \int_\GO \Big[\bpsi^\Gs_{\tau,m}\big(\bv_{\tau,m}+\bu_{\tau,m}-v-u\big)+\Big(g(x,u)+m\big|v-f(\by_{\tau,m})\big|^2+\tau\big|v-f(\by)\big|^2\\
\nnb &&\qq +\tau |u-\bu|^2+\Gs|v-\bv_{\tau,m}|^2+\Gs |u-\bu_{\tau,m}|^2\Big)\\
\nnb &&-\Big(g(x,\bu_{\tau,m})+m\big|\bv_{\tau,m}-f(\by_{\tau,m})\big|^2+\tau\big|\bv_{\tau,m}-f(\by)\big|^2+\tau|\bu_{\tau,m}-\bu|^2\Big)\Big]\, dx\\
&& \qq\qq\qq\qq\qq\all (v,u)\in \cV\times \cU.
\end{eqnarray}
Now, denote
$$
\eta_\Gve:=\max_{x\in \set{\na\by_{\tau,m}=0}}\big|\na\by^{\Gs,\Gve}_{\tau,m}(x)\big|.
$$
Then it follows from \refeq{E408} that $\ds\lim_{\Gve\to 0^+}\eta_\Gve=0$. On the other hand, by \refeq{E401A} and that  $\bpsi^{\Gs,\Gve}_{\tau,m}$ is bounded uniformly in $H^1_0(\GO)$ respect to $\Gs>0$ and $\Gve\in (0,1)$, we have
\begin{eqnarray*}
\nnb && (p-1)(\Gve^2+\eta_\Gve^2)^{p-2\over 2}\int_ {\set{\na\by_{\tau,m}=0}}\big|\nabla\bpsi^{\Gs,\Gve}_{\tau,m}\big|^2\, dx\\
\nnb &\leq & (p-1)\int_\GO  \Big(\Gve^2+|\nabla\by^{\Gs,\Gve}_{\tau,m}|^2\Big)^{p-2\over 2}\big|\nabla\bpsi^{\Gs,\Gve}_{\tau,m}\big|^2\, dx\\
\nnb &\leq & C_{\tau,m},
\end{eqnarray*}
where  $C_{\tau,m}$ is a constant independent of $\Gve\in (0,1)$ and $\Gs>0$.
Thus
\begin{equation}\label{E422}
\int_ {\set{\na\by_{\tau,m}=0}}\big|\nabla\bpsi^\Gs_{\tau,m}\big|^2\, dx\leq \lim_{\Gve\to 0^+} \int_ {\set{\na\by_{\tau,m}=0}}\big|\nabla\bpsi^{\Gs,\Gve}_{\tau,m}\big|^2\, dx=0.
\end{equation}
That is
\begin{equation}\label{E423}
 \na \bpsi^\Gs_{\tau,m}=0\qq\qq\qq\qq\qq\eqae \set{\na\by_{\tau,m}=0}.
\end{equation}

Finally, by \refeq{E419}, $\bpsi^\Gs_{\tau,m}$ is bounded uniformly in $H^1_0(\GO)$ respect to $\Gs>0$.
Then we can suppose that
$\bpsi^\Gs_{\tau,m}$ converges to some $\bpsi_{\tau,m}$ weakly in $H^1_0(\GO)$ and strongly in $L^2(\GO)$. By discussions similar to the above, we get
\refeq{E412C},\refeq{E413A}, \refeq{E414A} and \refeq{E415A} from  \refeq{E419}, \refeq{E420}, \refeq{E423} and \refeq{E421}, respectively.
\endpf

\section{Proof of the main result}
\def\theequation{5.\arabic{equation}}
\setcounter{equation}{0}

Similar to Proposition \ref{P402}, we have

\begin{Proposition}\label{P501}  Assume that {\em (S1)--(S3)} hold. Let
$(\bv_{\tau,m},\bu_{\tau,m})$ be an optimal control of Problem
$(P_{\tau,m})$ and $\by_{\tau,m}$ be the corresponding optimal state. Then, it holds that
as $m\to +\infty$,
\begin{equation}\label{E501}
\left\{\begin{array}{ll}
\ds \by_{\tau,m}\to \by & \ui{\,C^1(\bGO)}, \\
\ds \bv_{\tau,m}\to f(\by) & \si{\,L^2(\GO)}, \\
\ds \bu_{\tau,m}\to \bu & \si{\,L^2(\GO)}.
\end{array}\right.
\end{equation}
\end{Proposition}
\Proof
Similar to the proof of Proposition \ref{P402}, it suffices to prove \refeq{E501} in the sense of subsequence. We can suppose
that (as $m\to +\infty$)
\begin{equation}\label{E502}
\left\{\begin{array}{ll}
\ds \by_{\tau,m}\to \tiy & \ui{\,C^1(\bGO)}, \\
\ds \bv_{\tau,m}\to \tiv  & \wi{\,L^2(\GO)}, \\
\ds \bu_{\tau,m}\to \tiu  & \wi{\,L^2(\GO)}
\end{array}\right.
\end{equation}
for some $(\tiv,\tiu)\in \cV\times \cU$. Then
\begin{equation}\label{E503}
\left\{
\begin{array}{ll}
-\Div(|\nabla \tiy|^{p-2}\nabla \tiy)=\tiv+\tiu
&\eqin \GO,\\[1mm]
\tiy=0&\eqon
\pa \GO.
\end{array}
\right.
\end{equation}
On the other hand, it is easy to see that $y^{(f(\by),\bu)}=\by$. Then, by the optimality of $(\bv_{\tau,m},\bu_{\tau,m})$, it holds that
\begin{equation}\label{E504}
J_{\tau,m}(\bv_{\tau,m},\bu_{\tau,m})\leq J_{\tau,m}(f(\by),\bu)= J(\by,\bu).
\end{equation}
Thus, for any $M>0$, we have
\begin{eqnarray}\label{E505}
\nnb && \int_\GO \Big(f^0(x,\tiy)+g(x,\tiu)+M\left|\tiv-f(\tiy)\right|^2\Big)\, dx\\
\nnb && +\tau\Limsup_{m\to +\infty}\int_\GO \Big(\big|\bv_{\tau,m}-f(\by)\big|^2+|\bu_{\tau,m}-\bu|^2\Big)\, dx\\
\nnb &\leq & \Liminf_{m\to +\infty} \int_\GO \Big(f^0(x,\by_{\tau,m})+g(x,\bu_{\tau,m})+M\left|\bv_{\tau,m}-f(\by_{\tau,m})\right|^2\Big)\, dx\\
\nnb && +\tau\Limsup_{m\to +\infty}\int_\GO \Big(\big|\bv_{\tau,m}-f(\by)\big|^2+|\bu_{\tau,m}-\bu|^2\Big)\, dx\\
\nnb &\leq & \Limsup_{m\to +\infty} \int_\GO\Big(f^0(x,\by_{\tau,m})+g(x,\bu_{\tau,m})+m\left|\bv_{\tau,m}-f(\by_{\tau,m})\right|^2\\
\nnb && \qq +\tau\big|\bv_{\tau,m}-f(\by)\big|^2+\tau |\bu_{\tau,m}-\bu|^2\Big)\, dx\\
&= & \Limsup_{m\to +\infty} J_{\tau,m}(\bv_{\tau,m},\bu_{\tau,m})\leq J(\by,\bu).
\end{eqnarray}
This implies $\ds \int_\GO\big|\tiv-f(\tiy)\big|^2\, dx=0$. That is, $\tiv=f(\tiy)$ $\eqae \GO$. Consequently,
$(\tiy,\tiu)\in \cA$. Then, \refeq{E505} becomes
$$
 J(\tiy,\tiu)+\tau\Limsup_{m\to +\infty}\int_\GO \Big(\big|\bv_{\tau,m}-f(\by)\big|^2+|\bu_{\tau,m}-\bu|^2\Big)\, dx
\leq  J(\by,\bu).
$$
Since $J(\by,\bu)\leq J(\tiy,\tiu)$, we get
$$
\lim_{m\to +\infty}\int_\GO \big|\bv_{\tau,m}-f(\by)\big|^2\, dx=\lim_{m\to +\infty}\int_\GO |\bu_{\tau,m}-\bu|^2 \, dx=0.
$$
Then \refeq{E501} follows. \endpf

Now, we give a proof of our main theorem.

\textbf{Proof of Theorem \ref{T201}.}
By Proposition \ref{P501}, there is an $N_\tau\geq 1$ such that
\begin{equation}\label{E506}
\|f(\by_{\tau,m})-f(\by)\|_{C^1(\bGO)}< {1\over 2}, \qq\all m\geq N_\tau.
\end{equation}
We suppose that $m\geq N_\tau$ in the following.

One can easily see that \refeq{E415A} is equivalent to that both of the following inequalities hold:
\begin{eqnarray}\label{E507}
\nnb 0&\leq & \int_\GO \Big[\bpsi_{\tau,m}\big(\bv_{\tau,m}-v\big)+\Big(m\big|v-f(\by_{\tau,m})\big|^2+\tau\big|v-f(\by)\big|^2\Big)\\
 && -\Big(m\big|\bv_{\tau,m}-f(\by_{\tau,m})\big|^2+\tau\big|\bv_{\tau,m}-f(\by)\big|^2\Big)\Big]\, dx, \q \all v\in  \cV
\end{eqnarray}
and
\begin{eqnarray}\label{E508}
\nnb 0&\leq & \int_\GO \Big[\bpsi_{\tau,m}\big(\bu_{\tau,m}-u\big)+\Big(g(x,u)+\tau |u-\bu|^2\Big)\\
&&-\Big(|g(x,\bu_{\tau,m})+\tau|\bu_{\tau,m}-\bu|^2\Big)\Big]\, dx, \qq\all u\in \cU.
\end{eqnarray}
It is well-known that \refeq{E507} implies
\begin{eqnarray}\label{E509}
\nnb && \bpsi_{\tau,m}\bv_{\tau,m}-m\big|\bv_{\tau,m}-f(\by_{\tau,m})\big|^2-\tau\big|\bv_{\tau,m}-f(\by)\big|^2\\
\nnb &=& \max_{f(\by)-1\leq v\leq f(\by)+1}\Big(\bpsi_{\tau,m}v-m\big|v-f(\by_{\tau,m})\big|^2-\tau\big|v-f(\by)\big|^2\Big) ,\\
&& \hspace{7cm} \eqae \GO.
\end{eqnarray}
Therefore
\begin{equation}\label{E510}
\bv_{\tau,m}=\left\{\begin{array}{ll} v_{\tau,m},  & \eqif v_{\tau,m}\in [f(\by)-1,f(\by)+1], \\
f(\by)+1, & \eqif v_{\tau,m}>f(\by)+1,\\
f(\by)-1, & \eqif v_{\tau,m}<f(\by)-1,
\end{array}\right.
\end{equation}
where
\begin{equation}\label{E511}
v_{\tau,m}={1\over m+\tau}\big({1\over 2} \bpsi_{\tau,m}+mf(\by_{\tau,m})+\tau f(\by)\big).
\end{equation}
By \refeq{E511},
\begin{equation}\label{E512}
\bpsi_{\tau,m}= 2m \big(\bv_{\tau,m}-f(\by_{\tau,m})\big)+2\tau\big(\bv_{\tau,m}-f(\by)\big)\q
\eqon \set{\bv_{\tau,m}=v_{\tau,m}}.
\end{equation}
By  \refeq{E510}, it holds that
\begin{equation}\label{E513}
|\bv_{\tau,m}-f(\by)|=1, \q |v_{\tau,m}-f(\by)|>1, \qq \eqae \set{\bv_{\tau,m}\ne v_{\tau,m}}.
\end{equation}
Therefore, using \refeq{E506}, we have
\begin{eqnarray}
\nnb && |\bpsi_{\tau,m}|\geq |\bpsi_{\tau,m}+2m(f(\by_{\tau,m})-f(\by))|-m\\
\nnb &= &2(m+\tau)|v_{\tau,m}-f(\by)|-m\geq m\\
\nnb &\geq &{m\over 2} \big(|\bv_{\tau,m}-f(\by)|+|f(\by)-f(\by_{\tau,m})|\big)\\
\label{E514}&\geq & {m\over 2} |\bv_{\tau,m}-f(\by_{\tau,m})|, \qq\eqae \set{\bv_{\tau,m}\ne v_{\tau,m}}.
\end{eqnarray}
Combining \refeq{E512} with \refeq{E514}, we get
\begin{equation}\label{E515}
\big|m \big(\bv_{\tau,m}-f(\by_{\tau,m})\big)\big|\leq
2|\bpsi_{\tau,m}|+\tau, \qq \eqae \GO.
\end{equation}
Then it follows easily from \refeq{E412C} and \refeq{E515} that
\begin{eqnarray}\label{E516}
\nnb &&\|\bpsi_{\tau,m}\|^2_{H^1_0(\GO)}+\int_{\set{\nabla\by_{\tau,m}\ne 0}}|\nabla\by_{\tau,m}|^{p-2}\, |\nabla\bpsi_{\tau,m}|^2\, dx\\
\nnb &\leq & C\Big(1+\Big|\int_\GO 2m \big(f(\by_{\tau,m})-\bv_{\tau,m}\big)f^\pri(\by_{\tau,m})\bpsi_{\tau,m}\, dx\Big|\Big)\\
&\leq & \tiC+\tiC\int_\GO|\bpsi_{\tau,m}|^2\,dx
\end{eqnarray}
for some constant $\tiC>0$ independent of $m\ge 1$ and $\tau\in (0,1)$.

Let
$$
\mu_{\tau,m}:={1\over \|\bar{\psi}_{\tau,m}\|_{L^2(\GO)}+1}, \q \Phi_{\tau,m}:=\mu_{\tau,m}\bpsi_{\tau,m}.
$$
Then
\begin{equation}\label{E517}
\mu_{\tau,m}\geq 0, \q \|\Phi_{\tau,m}\|_{L^2(\GO)}+\mu_{\tau,m}=1
\end{equation}
and
\begin{equation}\label{E518}
\|\Phi_{\tau,m}\|^2_{H^1_0(\GO)}+\int_{\set{\nabla\by_{\tau,m}\ne 0}}|\nabla\by_{\tau,m}|^{p-2}\, |\nabla\Phi_{\tau,m}|^2\, dx\leq 2\tiC.
\end{equation}
Thus, along a subsequence of $m\to +\infty$, we have $\mu_{\tau,m}\to \mu_\tau$ and
\begin{equation}\label{E519}
\Phi_{\tau,m}\to \bpsi_\tau \qq\wi{H^1_0(\GO)}, \si{L^2(\GO)}
\end{equation}
for some constant $\mu_\tau$ and $\bpsi_\tau\in H^1_0(\GO)$.  By \refeq{E517} and \refeq{E518},  we have
\begin{equation}\label{E520}
\mu_\tau\geq 0, \q \|\bpsi_\tau\|_{L^2(\GO)}+\mu_\tau=1
\end{equation}
and
\begin{equation}\label{E521}
\|\bpsi_\tau\|_{H^1_0(\GO)}^2+\int_{\set{\nabla\by\ne 0}}|\nabla\by|^{p-2}\, |\nabla\bpsi_\tau|^2\, dx\leq  2\tiC.
\end{equation}

On the other hand, by \refeq{E518} and Sobolev's inequality, $\ds \|\Phi_{\tau,m}\|_{L^q(\GO)}$ is bounded uniformly for some $q>2$. While by \refeq{E501} and
\refeq{E513},  the  Lebesgue measure of $\set{\bv_{\tau,m}\ne v_{\tau,m}}$ (denote it by $\big|{\set{\bv_{\tau,m}\ne v_{\tau,m}}}\big|$) tends to zero as $m\to +\infty$. Thus,  as $m\to +\infty$,
\begin{eqnarray}\label{E522}
\nnb && \|2\mu_{\tau,m}m \big(\bv_{\tau,m}-f(\by_{\tau,m})\big)\chi_{\set{\bv_{\tau,m}\ne v_{\tau,m}}}\|_{L^2(\GO)}\\
\nnb &\leq & \Big\|\Big(4\big|\Phi_{\tau,m}\big|+2\Big)\chi_{\set{\bv_{\tau,m}\ne v_{\tau,m}}}\Big\|_{L^2(\GO)}\\
&\leq & \big|{\set{\bv_{\tau,m}\ne v_{\tau,m}}}\big|^{q-2\over q}\, \Big\|\Big(4\big|\Phi_{\tau,m}\big|+2\Big)\Big\|_{L^q(\GO)}\to 0.
\end{eqnarray}
Combining the above with \refeq{E512}, we get that as $m\to +\infty$,
\begin{equation}\label{E523}
2\mu_{\tau,m}m \big(\bv_{\tau,m}-f(\by_{\tau,m})\big)\to \bpsi_\tau , \si {L^2(\GO)}.
\end{equation}
Thus, it follows from \refeq{E501},  \refeq{E413A}--\refeq{E414A}, \refeq{E508} and \refeq{E523} that
\begin{eqnarray}\label{E524}
\nnb && \ds
-\Div\Big[|\nabla\by|^{p-2}\Big(I+(p-2)
\frac{\nabla\by(\nabla\by)^T}
{|\nabla\by|^2}\Big)\nabla\bpsi_\tau\Big]\\
 &=& -\mu_\tau f^0_y(x,\by)
+f^\pri(\by)\bpsi_\tau, \qq\q\qq \eqin \set{\na\by\ne 0},\\
\nnb &&\\
\label{E525} && \na \bpsi_\tau=0\qq\qq\qq\qq\qq\eqae \set{\na\by=0}
\end{eqnarray}
and
\begin{eqnarray}\label{E526}
 0&\leq & \int_\GO \Big[\bpsi_\tau\big(\bu-u\big)+\mu_\tau \big(g(x,u)-g(x,\bu)+\tau |u-\bu|^2 \big)\Big]\, dx, \, \all u\in \cU.
\end{eqnarray}
Finally, by \refeq{E520}---\refeq{E521},  we can suppose that as $\tau\to 0^+$, it holds that  $\mu_\tau\to \mu$ and
\begin{equation}\label{E527}
\bpsi_\tau\to \bpsi \qq\wi{H^1_0(\GO)}, \si{L^2(\GO)}
\end{equation}
for some constant $\mu$ and  $\bpsi\in H^1_0(\GO)$. Then it follows from \refeq{E521}
and \refeq{E524}---\refeq{E527} that
\begin{equation}\label{E528}
\mu\geq 0, \q \|\bpsi\|_{L^2(\GO)}+\mu=1,
\end{equation}
\refeq{E204}---\refeq{E206} and $\bpsi\in \cX$.
 We get the proof.
\endpf

\begin{Remark}
From the proof of Theorem \ref{T201}, we can see that if $\ds \|\bpsi_{\tau,m}\|_{L^2(\GO)}$ is bounded uniformly, then $\mu>0$. In particular, if $f^\pri(y)\leq 0$, which is the case that the state equation \refeq{E201} is well-defined,  then $\mu>0$.
\end{Remark}

\section{An Example}
\def\theequation{6.\arabic{equation}}
\setcounter{equation}{0}
We give a simple example to show a usage of our
main theorem.
\begin{Example}
Assume that:
\begin{enumerate}[{\rm (i)}]
\item the function $f$ satisfies $(S2)$ and
\begin{equation}\label{E601}
f(y)+a>0, \qq\all y\in \IR;
\end{equation}
\item the function $f^0(x,y)\equiv f^0(y)$ satisfies $(S3)$;
\item the set $\ds \set{y\in\IR\big| f^\pri(y)=f^0_y(y)}$ is at most countable;
\item $g(x,u)=u$.
\end{enumerate}
Let $(\by,\bu)$ be an optimal pair to Problem $(P)$. Then $\bu$ should be a bang-bang control.
\end{Example}
\proof
In this example, \refeq{E206} becomes
 \begin{equation}\label{E606}
(\bpsi(x) -\mu) \bu(x)=\max_{a\leq u\leq b}(\bpsi(x) -\mu)  u , \qq\eqae x\in \GO.
\end{equation}
Thus, to prove $\bu$ is bang-bang, it need only to prove $\set{\bpsi=\mu}$ has zero measure.

By Assumption (iii), we can suppose that
\begin{equation}\label{E603}
\set{y\in\IR\big| f^\pri(y)=f^0_y(y)}\subseteq \set{y_1,y_2,\ldots}.
\end{equation}
It is well-known that (see \cite{Mo1}) if $y\in W^{1,1}_{loc}(\GO)$, then for any constant $C$,
\begin{equation}\label{E604}
\na y(x)=0, \qq \eqae x\in \set{y=C}.
\end{equation}
Here, it is crucial that there holds a similar
. By Theorem 1.1 in
 \cite{Lou6},
\begin{equation}\label{E605}
\na \by(x)=0, \qq \eqae x\in \set{\na \by=0}.
\end{equation}
Consequently, combining the above with \refeq{E601}, we can see that  $\set{\na\by=0}$ has zero measure.

Now, if  $\set{\bpsi=\mu}$ has positive measure, then by
\refeq{E604},
$$
\na \bpsi(x)=0, \qq \eqae x\in \set{\bpsi=\mu}.
$$
Noting
that $\by\in C^{1,\Ga}(\bGO)$, we have $\bpsi\in W^{2,2}_{loc}\big(\set{\na \by\ne 0}\big)$. Therefore,
\begin{eqnarray}\label{E607}
\nnb & & \mu\big(f^\pri(\by(x))-f^0_y(\by(x))\big)=f^\pri(\by(x))\bpsi(x)-\mu f^0_y(\by(x))\\
 &=& -\Div \Big[|\nabla\by|^{p-2}\Big(I+(p-2)
\frac{\nabla\by(\nabla\by)^T}{|\nabla\by|^2}\Big)\na \bpsi\Big]=0 , \qq\eqae x\in \set{\bpsi=\mu}.
\end{eqnarray}
By \refeq{E604},
$$
\na \by(x)=0, \qq \eqae x\in \set{\by=y_k}.
$$
Thus  $\ds\set{\by=y_k}$ has zero measure and consequently $\ds \set{f^\pri(\by)=f^0_y(\by)}$ has zero measure.
Therefore, by \refeq{E607} and that $\ds \set{\bpsi=\mu}$ has positive measure, we get $\mu=0$ and $\ds\set{\bpsi=0}$ has positive measure.  Moreover,
\refeq{E204} becomes
\begin{eqnarray}\label{E608}
& & -\Div \Big[|\nabla\by|^{p-2}\Big(I+(p-2)
\frac{\nabla\by(\nabla\by)^T}{|\nabla\by|^2}\Big)\na \bpsi\Big]= f^\pri(\by)\bpsi, \qq\eqin \{\na \by\neq 0\}.
\end{eqnarray}
Since $\ds\set{\bpsi=0}$ has positive measure, by Proposition 3 in \cite{Figueiredo-Gossez}, $\bpsi$ has a zero $x_0$ of infinite order, i.e., for any $m$,
$$
\int_{|x-x_0|\leq r}|\bpsi(x)|^2\, dx=O(r^m), \qq r\to 0^+.
$$
Then, by Theorem 1.1 in \cite{Garofalo-Lin}, we get that
$$
\bpsi(x)=0, \qq\eqae x\in\GO.
$$
Contradicts to \refeq{E203}. Therefore, $\ds\set{\bpsi=\mu}$ has zero measure and $\bu$ is a bang-bang control.
\endpf

\end{document}